\documentclass{article}
\usepackage[leqno]{amsmath}
\usepackage{mathrsfs}
\usepackage{graphics,amssymb,amsthm,amsmath}   
\usepackage{latexsym}   
\usepackage{eufrak}
\usepackage{epsfig}
\usepackage{euscript}
\usepackage{color}
\usepackage{url}
\usepackage{graphicx,psfrag}
\usepackage{endnotes}

\allowdisplaybreaks

\theoremstyle{remark}
\newtheorem{remark}{\bf Remark}[section]

\newtheorem{discussion}{\bf Discussion}[section]

\theoremstyle{definition}
\newtheorem{defn}[remark]{Definition}

\theoremstyle{theorem}
\newtheorem{thm}[remark]{Theorem}

\definecolor{thmcolor}{rgb}{0,0,.4}
\definecolor{remarkcolor}{rgb}{0,.2,0} 
\definecolor{proofcolor}{rgb}{.4,0,0} 
\definecolor{quecolor}{rgb}{.2,.2,0} 
\definecolor{axcolor}{rgb}{.23,0,.23} 
\definecolor{red}{rgb}{1,0,0} 

\newcommand{\comment}[1]{}
\newcommand{\rmk}[1]{} 

\newcommand{\de}{\ \overset{d}{=}\ } 

\newcommand{\SCA}{\mbox{\sf SCA\/}} %
\newcommand{\CA}{\mbox{\sf CA\/}} %
\newcommand{\Df}{\mbox{\sf Df\/}} %
\newcommand{\RDf}{\mbox{\sf RDf\/}} %
\newcommand{\RA}{\mbox{\sf RA\/}} %
\newcommand{\SA}{\mbox{\sf SA\/}} %
\newcommand{\QRA}{\mbox{\sf QRA\/}} %
\renewcommand{\c}{\mbox{\sf c\/}} %
\newcommand{\s}{\mbox{\sf s\/}} %
\newcommand{\p}{\mbox{\sf p\/}} %
\newcommand{\q}{\mbox{\sf q\/}} %
\renewcommand{\d}{\mbox{\sf d\/}} %
\newcommand{\ZF}{\mbox{\it ZF\/}} %
\renewcommand{\L}{\mathcal{L}} %
\newcommand{\Ax}{\mbox{\sf Ax\/}} %
\newcommand{\SAx}{\mbox{\sf SAx\/}} %
\newcommand{\pair}{\mbox{\sf pair\/}} %
\newcommand{\comp}{\mbox{{\small{$\odot$}}}} %
\def\udot{{\cup\hspace*{-1.43mm} {\boldsymbol \cdot} }}
\newcommand{\conv}{^{\udot}} %
\newcommand{\inv}{^{\smallsmile}} %
\newcommand{\id}{\mbox{\tiny{$\overset{{\boldsymbol \cdot} }{1}$\/}}\/} %
\newcommand{\mi}{-} %
\newcommand{\dra}{\mbox{\sf Dra\/}} %
\newcommand{\Dra}{{\mathfrak Dra\/}} %
\newcommand{\dca}{\mbox{\sf Dca\/}} %
\newcommand{\Dca}{{\mathfrak Dca\/}} %
\newcommand{\Tr}{{\sf Tr\/}} %
\newcommand{\tr}{{\sf tr\/}} %
\newcommand{\h}{{\sf h\/}} %
\newcommand{\ZFC}{\mbox{\it ZF\/}} %
\newcommand{\op}{\mbox{\sf op\/}} %
\newcommand{\Triplet}{{\sf Triplet\/}} %
\newcommand{\T}{{\sf T\/}} %
\newcommand{\D}{{\sf d\/}} %
\newcommand{\e}{{\sf e\/}} %
\renewcommand{\t}{{\sf t\/}} %
\newcommand{\Fm}{\mbox{{\it Fm\/}}} %
\newcommand{\inn}{{\varepsilon\/}} %
\newcommand{\In}{\epsilon} 
\renewcommand{\)}{\mbox{)\!\!)\/}}  
\newcommand{\bz}{\mbox{)\!)\/}}  
\newcommand{\nz}{\mbox{(\!(\/}} 
\newcommand{\Fmd}{\mbox{\it Fmd\/}}
\newcommand{\Ld}{\mbox{$\mathcal Ld$\/}}
\newcommand{\f}{\varphi}
\newcommand{\gFm}{\mathfrak {Fm}} 
\newcommand{\gFmd}{\mathfrak {Fmd}} 
\newcommand{\Mm}{\mathfrak M}
\newcommand{\Cc}{\mathfrak C}
\newcommand{\gFr}{\mathfrak {Fr}} 
\newcommand{\gZd}{\mathfrak {Zd}} 
\newcommand{\Val}{\mbox{\sf Val\/}} %
\newcommand{\val}{\mbox{\sf val\/}} %
\newcommand{\lr}{\leftrightarrow}
\newcommand{\esec}{\mbox{$\langle\rangle$}} 

\newcommand{\True}{\mbox{\sf True\/}} %
\newcommand{\exy}{\mbox{$x\!\doteq\! y$\/}} %
\newcommand{\exz}{\mbox{$x\!\doteq\! z$\/}} %
\newcommand{\eyz}{\mbox{$y\!\doteq\! z$\/}} %
\newcommand{\eyx}{\mbox{$y\!\doteq\! x$\/}} %
\newcommand{\ezx}{\mbox{$z\!\doteq\! x$\/}} %
\newcommand{\ezy}{\mbox{$z\!\doteq\! y$\/}} %
\newcommand{\exx}{\mbox{$x\!\doteq\! x$\/}} %
\newcommand{\eyy}{\mbox{$y\!\doteq\! y$\/}} %
\newcommand{\ezz}{\mbox{$z\!\doteq\! z$\/}} %
\newcommand{\euv}{\mbox{$u\!\doteq\! v$\/}} %
\newcommand{\eq}{\!\doteq\!}

\newcommand\undharom{\,|\hspace*{-3.5pt}\begin{array}{@{\hspace*{2pt}}c@{\hspace*{2pt}}}
\phantom{\scriptstyle r,\alpha}\\[-2.5mm]
\text{---}\!\!\text{---}\\[-3mm]
\scriptstyle d
\end{array}}

\newcommand{\vdr}{\dr} 
\newcommand{\vd}{{\tiny{\undharom}}}
\newcommand{\undra}{\,|\hspace*{-3.5pt}\begin{array}{@{\hspace*{2pt}}c@{\hspace*{2pt}}}
\phantom{\scriptstyle r,\alpha}\\[-1mm]
\text{---}\!\!\text{---}\\[-1mm]
\scriptstyle {\ 3}
\end{array}}
\newcommand{\undre}{\,|\hspace*{-3.5pt}\begin{array}{@{\hspace*{2pt}}c@{\hspace*{2pt}}}
\phantom{\scriptstyle r,\alpha}\\[-1mm]
\text{---}\!\!\text{---}\\[-1mm]
\scriptstyle {\ r}
\end{array}}
\newcommand{\undrn}{\,|\hspace*{-3.5pt}\begin{array}{@{\hspace*{2pt}}c@{\hspace*{2pt}}}
\phantom{\scriptstyle r,\alpha}\\[-1mm]
\text{---}\!\!\text{---}\\[-1mm]
\scriptstyle {\ n}
\end{array}}
\newcommand{\undrf}{\,|\hspace*{-3.5pt}\begin{array}{@{\hspace*{2pt}}c@{\hspace*{2pt}}}
\phantom{\scriptstyle r,\alpha}\\[-1mm]
\text{---}\!\!\text{---}\\[-1mm]
\scriptstyle {\ 4}
\end{array}}
\newcommand{\dr}{\tiny{\undra}} 
\newcommand{\vr}{\tiny{\undre}}
\newcommand{\vn}{\tiny{\undrn}}
\newcommand{\vf}{\tiny{\undrf}}

\newcommand{\notvd}{\mbox{$\quad\not\!\!\!\!\!\!\vd$}}
\newcommand{\notmodels}{\mbox{$\,\not\!\models\,$}}

\begin{document}

\title{Formalizing set theory in weak logics, searching for the
weakest logic with G\"odel's incompleteness property.}
\author{Hajnal Andr\'eka and Istv\'an N\'emeti
\thanks{Research supported by the
        Hungarian National Foundation for scientific research grants
        T81188.}}
\date{2011. July}
\maketitle

\begin{abstract}
We show that first-order logic can be translated into a very simple
and weak logic, and thus set theory can be formalized in this weak
logic. This weak logical system is equivalent to the equational
theory of Boolean algebras with three commuting complemented closure
operators, i.e., that of diagonal-free 3-dimensional cylindric
algebras ($\Df_3$'s). Equivalently, set theory can be formulated in
propositional logic with 3 commuting S5 modalities (i.e., in the
multimodal logic [S5,S5,S5]). There are many consequences, e.g.,
free finitely generated $\Df_3$'s are not atomic and [S5,S5,S5] has
G\"odel's incompleteness property. The results reported here are
strong improvements of the main result of the book: Tarski, A. and
Givant, S. R., Formalizing Set Theory without variables, AMS, 1987.
\end{abstract}

\section{Introduction}\label{intro-sec}

Tarski in 1953 \cite{T53a, T53b} formalized set theory in the theory
of relation algebras. Why did he do this? Because the equational
theory of relation algebras ($\RA$) corresponds to a logic without
individual variables, in other words, to a propositional logic. This
is why the title of the book \cite{TG} is ``Formalizing set theory
without variables". Tarski got the surprising result that a
propositional logic can be strong enough to ``express all of
mathematics", to be the arena for mathematics. The classical view
before this result was that propositional logics in general were
weak in expressive power, decidable, uninteresting in a sense. By
using the fact that set theory can be built up in it, Tarski proved
that the equational theory of $\RA$ is undecidable. This was the
first propositional logic shown to be undecidable.

From the above it is clear that replacing $\RA$ in Tarski's result
with a ``weaker" class of algebras is an improvement of the result
and it is worth doing. For more on this see
Tarski-Givant~\cite[pp.$89_2-90^4$ and footnote 17 on p.90]{TG},
especially the open problem formulated therein.

A result of J.\ D.\ Monk says that for every finite $n$ there is a
3-variable first-order logic (FOL) formula which is valid but which
can be proved (in FOL) with more than $n$ variables only.
Intuitively this means that during any proof of this formula there
are steps when we have to use $n$ independent data (stored in the
$n$ variables as in $n$ machine
registers). 
For example, the associativity of relation composition of binary
relations can be expressed with 3 variables but 4 variables are
needed for any of its proofs.

Tarski's main idea in \cite{TG} is to use pairing functions to form
ordered pairs, and so to store two pieces of data in one register.
He used this technique to translate usual infinite-variable
first-order logic FOL into the three-variable fragment of it. From
then on, he used that any three-variable FOL-formula about binary
relations can be expressed by an $\RA$-equation, \cite[sec
5.3]{HMTII}. He needed two registers for storing the data belonging
to a binary relation and he had one more register available for
making computations belonging to a proof.

The finite-variable fragment hierarchy of FOL corresponds to
cylindric algebras (\CA's). The $n$-variable fragment $\L_n$ of FOL
consists of all FOL-formulas which use only the first $n$ variables.
By Monk's result, $\L_n$ is essentially incomplete for all $n\ge 3$,
it cannot have a finite Hilbert-style complete and sound inference
system. We get a finite Hilbert style inference system $\vn$ for
$\L_n$ by restricting a usual complete one for infinite-variable FOL
to the first $n$ variables (see \cite[sec. 4.3]{HMTII}). This
inference system $\vn$ belonging to $\L_n$ expresses $\CA_n$, it is
sound but not complete: $\vn$ is much weaker than validity
$\models_n$.

Relation algebras are halfway between $\CA_3$ and $\CA_4$, the
classes of 3-dimen\-sional and 4-dimensional cylindric algebras,
respectively. We sometimes jokingly say that $\RA$ is $\CA_{3.8}$.
Why is $\RA$ stronger than $\CA_3$? Because, the so-called relation
algebra reduct of a $\CA_3$ is not necessarily an $\RA$, e.g.,
associativity of relation composition can fail in the reduct. See
\cite[sec 5.3]{HMTII}, and for more in this line see
N\'emeti-Simon~\cite{NeSiIGPL97}. \rmk{The class $\SA$, defined by
weakening the associativity, corresponds closely to $\CA_3$. Write
more, about Maddux's work.??} Why is $\CA_4$ stronger than $\RA$?
Because not every $\RA$ can be obtained, up to isomorphism, as the
relation algebra reduct of a $\CA_4$
. However, the same equations are true in $\RA$ and in the
class of all relation algebra reducts of $\CA_4$'s (
Maddux's result, see \cite[sec 5.3]{HMTII}). Thus Tarski formulated
Set Theory, roughly, in $\CA_4$, i.e., in $\L_4$ with $\vf$, or in
$\L_3$ with validity $\models$.

N\'emeti \cite{NPrep}, \cite{NDis} improved this result by
formalizing set theory in $\CA_3$, i.e., $\L_3$ with $\vdr$ in place
of validity $\models$.
The main idea for this improvement was using the paring functions to
store all data always, during every step of a proof, in one register
only and so one got two registers to work with in the proofs. In
this approach one represents binary relations as unary ones (of
pairs). For the ``execution" of this idea see
sections~\ref{qra-sec}-\ref{set-sec} of the present paper.

First-order logic has equality as a built-in relation. One of the
uses of equality in FOL is that it can be used to express (simulate)
substitutions of variables, thus to ``transfer" content of one
variable to the other. The reduct $\SCA_3$ of $\CA_3$ ``forgets"
equality $\d_{ij}$ but retains substitution in the form of the
term-definable operations $\s^i_j$. The logic belonging to $\SCA_3$
is weaker than 3-variable fragment of FOL. Zal\'an Gyenis
\cite{Gyenis} improved parts of N\'emeti's result by using $\SCA_3$
in place of $\CA_3$.

We get a much weaker logic by forgetting substitutions, too, this is
the logic corresponding to $\Df_3$ in which we formalize Set Theory
in the present paper. Without equality or substitutions, if one has
only binary relations, one cannot really use the third variable for
anything; and it is known that the two-variable fragment of FOL is
decidable, so it is already too weak for formalizing set theory.
Therefore we need at least one ternary relation symbol (or atomic
formula) in order to use the third variable, while in the language
of set theory we only have one binary relation symbol, the
elementhood-relation $\In$. Therefore, while in formalizing set
theory in the three-variable fragment of FOL (in $\CA_3$) we could
do with one binary relation symbol, we did not have to change
vocabulary during the formalization, in the present equality- and
substitution-free case we have to change vocabulary, and we have to
pay attention to this new feature of the translation mapping. A key
device of our proofs will be a recursive ``translation mapping"
translating FOL into the equational language of $\Df_3$, or
equivalently into the logic $\Ld_3$ defined in section~\ref{fmd-sec}
below.

$\Df_3$ is nothing more than Boolean algebras with three commuting
complemented closure operators. The only connection between these
operators is commutativity. We know that without commutativity the
class is too weak for supporting set theory because its equational
theory is decidable \cite{NDis}. We know that two commuting such
operators do not suffice, for the same reason. We do not know how
much complemented-ness of the closure operators is important for
supporting set theory.

In section \ref{fmd-sec} we introduce our simple logic $\Ld_3$ in
several different forms, which reveal its propositional logic
character. Then we state three of the main theorems about this
logic: it is only seemingly weak, because set theory can be built up
in it (Thm.\ref{zfc-t}), and also G\"odel's incompleteness theorem
holds for it (Thm.\ref{incomp-t}). In contrast with the fact that
$\Ld_3$ cannot have a sound and complete Hilbert-style proof system,
we state a completeness theorem for $\Ld_3$ which comes very close
to having a Hilbert-style sound and complete proof system
(Thm.\ref{comp-t}). Sections \ref{qra-sec}-\ref{set-sec} contain a
full proof for Thm.\ref{zfc-t}, and a proof for a weaker version of
Thm.\ref{comp-t}. In section \ref{free-sec} we prove, as a corollary
of Thm.\ref{zfc-t}, that the finitely generated $\Df_3$'s are not
atomic. This proof also contains the main ideas for a proof of
Thm.\ref{incomp-t}.

We make the paper available in the present form because so many
people expressed strong interest in the proofs of two of the main
theorems, Thm.\ref{zfc-t} and Thm.\ref{free-t}. We will keep
developing the paper and new versions will be found on our
home-page, via the link
\url{http://www.renyi.hu/~nemeti/FormalizingST.htm}. Via that link
one can find more on the history of the problem settled in the
present paper, see \cite{shortnote}, and some unpublished works, see
\cite{NDis}, \cite{NPrep}.


\section{\label{fmd-sec} A simple logical system: three-variable logic without
equality or substitutions} In this section we define the ``target
logic" $\Ld_3$ of our translation. We give several different forms
for it to give a feeling of its expressive power. After this, we
formulate three of our main theorems, all stating unexpected
properties of this logic.

The language of our system contains three variable symbols, $x,y,z$,
one ternary relational symbol $P$, and only one atomic formula,
namely $P(x,y,z)$. (We note that, e.g., the formula $P(y,x,z)$ is
not available in this language.) The logical connectives are $\lor,
\neg, \exists x, \exists y, \exists z$. We denote the set of
formulas (of $\Ld_3$) by $\Fmd_3$. We will use the derived
connectives $\forall, \land, \to, \leftrightarrow$, too, as
abbreviations: $\forall v\f\de\neg\exists v\neg\f$,
$\f\land\psi\de\neg(\neg\f\lor\neg\psi)$, $\f\to\psi\de
\neg\f\lor\psi$, $\f\leftrightarrow\psi\de
(\f\to\psi)\land(\psi\to\f)$. Sometimes we will write, e.g.,
$\exists xy$ or $\forall xyz$ in place of $\exists x\exists y$ or
$\forall x\forall y\forall z$, respectively. $\Fmd^1_3$ denotes the
set of formulas in $\Fmd_3$ with one free variable $x$, we will
often deal with these in section~\ref{qra-sec} on.

The proof system $\vd$ which we will use is a Hilbert style one with
the following logical axioms and rules.

The logical axioms are the following.
Let $\varphi, \psi \in \Fmd_3$ and $v,w \in \{ x,y,z\}$.\\
((1))\quad $\varphi$, if $\varphi$ is a propositional tautology.\\
((2))\quad $\forall  v (\varphi \to \psi) \to (\exists  v
\varphi \to \exists  v \psi)$.\\
((3))\quad $\varphi \to \exists v \varphi$.\\
((4))\quad $\exists v \exists v \varphi \to \exists v\varphi$.\\
((5))\quad $\exists v (\varphi\lor\psi)\leftrightarrow (\exists v\varphi\lor\exists v\psi)$.\\
((6))\quad $\exists v\neg\exists v\varphi\to \neg\exists v\varphi$.\\
((7))\quad $\exists v\exists w\varphi\to\exists w\exists v\varphi$.\\
%
The inference rules are Modus Ponens ((MP), or detachment), and
Generalization~((G)).

This proof system is a direct translation of the equational axiom
system of $\Df_3$. Axiom ((2)) is needed for ensuring that the
equivalence relation defined on the formula algebra by $\f\equiv\psi
\Leftrightarrow \vd\f\!\leftrightarrow\!\psi$ be a congruence with
respect to (w.r.t.) the operation $\exists v$. It is congruence
w.r.t.\ the Boolean connectives $\lor,\lnot$ by axiom ((1)). Axiom
((1)) expresses that the formula algebra factorized with $\equiv$ is
a Boolean algebra, axiom ((5)) expresses that the quantifiers
$\exists v$ are operators on this Boolean algebra (i.e., they
distribute over $\lor$), axioms ((3)),((4)) express that these
quantifiers are closure operations, axiom ((6)) expresses that they
are complemented closure operators (i.e., the negation of a closed
element is closed again). Together with ((5)) they imply that the
closed elements form a Boolean subalgebra, and hence the quantifiers
are normal operators (i.e., the Boolean zero is a closed element).
Finally, axiom ((7)) expresses that the quantifiers commute with
each other.

We define $\Ld_3$ as the logic with formulas $\Fmd_3$ and with proof
system $\vd$. The logic $\Ld_3$ inherits a natural semantics from
first-order logic (FOL). The proof system $\vd$ is sound with
respect to this semantics, but it is not complete. Moreover, there
is no finite Hilbert-style inference system which would be complete
and sound at the same time w.r.t.\ this semantics (because the
quasi-equational theory of $\RDf_3$ is not finitely axiomatizable,
see \cite{HMTII} and \cite{HbPhL}).

We note that in the above system, axiom ((6)) can be replaced with
the following ((8)):\\
((8))\quad $\exists v(\varphi\land\exists v\psi)\leftrightarrow (\exists v\varphi\land\exists v\psi)$.\\

\bigskip

In the present paper we will use our logic $\Ld_3$ as introduced
above. However, it has several different but equivalent forms, each
of which has advantages and disadvantages. We review some of the different forms below.

Restricted 3-variable FOL is introduced in \cite[Part II,
p.157]{HMTII}, with proof system $\vr$. If we restrict this system
$\vr$ to formulas not containing the equality $=$ then we get a
system equivalent to our $\Ld_3$. Lets call this system {\it
restricted 3-variable FOL without equality}.
That is, the formulas are those of restricted 3-variable FOL which
contain no equality, and we leave out from the axioms of $\vr$ the
axioms which contain equality. This way we get a proof system with
Modus Ponens and
Generalization as deduction rules and with the following axioms:\\

\noindent
((V1))\quad $\varphi$, if $\varphi$ is a propositional tautology.\\
((V2))\quad $\forall  v (\varphi \to \psi) \to (\forall  v
\varphi \to \forall  v \psi)$.\\
((V3))\quad $\forall  v\varphi \to \varphi$.\\
((V4))\quad $\varphi \to \forall  v \varphi$, if $v$ does not
occur free in~$\varphi$.\\

Lets call this\footnote{We note that we also omitted ((V9)) of
\cite{HMTII} because $\forall v$ abbreviates $\lnot\exists v\lnot$
in our approach, so ((V9)) is not needed.} {\it equality-free
$\vr$}. Rule ((V4)) in this system essentially uses individual
variables in its using the notion of free variables of a formula. On
the other hand, no axiom in $\vd$ needs to use the structure of a
formula occurring in a rule, it is essentially variable-free. So, an
advantage of $\vd$ over equality-free $\vr$ is that it is more
``algebraic", more like propositional logic. On the other hand,
equality-free $\vr$ contains fewer axioms (it contains only
((V1))-((V4)) as axioms).

The logic $\Ld_3$ has a neat {\it modal logic form}: three commuting
S5 modalities. This is denoted as [S5,S5,S5], see \cite[p.379, lines
15-20]{GKWZ}. We recall this logic in a slightly simplified form.
The language contains one propositional variable $p$, the
connectives are $\lor, \neg, \Diamond_1, \Diamond_2, \Diamond_3$. We
use $\Box_i\de\neg\Diamond_i\neg$, $\to, \leftrightarrow$ as derived
connectives as before, and the axioms are the following (where $\f,
\psi$ are arbitrary formulas of the language and $i,j\in\{ 1,2,3\}$):\\

\noindent
((B))\quad $\varphi$, if $\varphi$ is a propositional tautology,\\
((K))\quad $\Box_i(\f\to\psi)\to(\Box_i\f\to\Box_i\psi)$,\\
((S5))\quad $\Diamond_i\f\to\Box_i\Diamond_i\f$,\\
((C1))\quad $\Diamond_i\Diamond_j\f\to\Diamond_j\Diamond_i\f$,\\
((C2))\quad $\Diamond_i\Box_j\f\to\Box_j\Diamond_i\f$.\\

\noindent The rules are Modus Ponens and Generalization (or, in
other word, Necessitation, i.e., $\f\vdash\Box_i\f$). This modal
logic is complete w.r.t.\ the frames consisting of three commuting
equivalence relations as accessibility relations for the three
modalities.

One can present this logic in yet one different form: {\it
Equational logic} as the background logic, and the defining axioms
of $\Df_3$ as logical axioms. (Occasionally, we refer to this logic
informally as ``the equational theory of $\Df_3$".) For
completeness, we include this form of $\Ld_3$ here, too. The
language consists of equations $\tau = \sigma$ where $\tau, \sigma$
are terms built up from (arbitrarily many) variables by the use of
the function symbols $+,-,f,g,h$ where $+$ is binary and the rest
are unary. The axioms
are the following, where $x,y,z$ are variables and $F\in\{ f,g,h\}$:\\

\noindent
((B1))\quad $x+y=y+x$, \\
((B2))\quad $x+(y+z) = (x+y)+z$,\\
((B3))\quad $-(-(x+y)+-(x+-y))=x$,\\
((D1))\quad $x+Fx=Fx$, \\
((D2))\quad $FFx=Fx$, \\
((D3))\quad $F(x+y)=Fx+Fy$, \\
((D4))\quad $F(-Fx)=-Fx$, \\
((D5))\quad $fgx=gfx$,\quad $fhx=hfx$,\quad $ghx=hgx$.\\

\noindent The rules are those of the equational logic:

\noindent Rules of equivalence:\\ $\tau=\tau$,\quad
$\tau=\sigma\vdash\sigma=\tau$,\quad $\tau=\sigma,
\sigma=\rho\vdash\tau=\rho$,\\
Rules of congruence:\\ $\tau=\sigma, \rho=\delta\quad\vdash\quad
-\tau=-\sigma, f\tau=f\sigma, g\tau=g\sigma, h\tau=h\sigma$,
$\tau+\rho=\sigma+\delta$,\\
Rule of invariance:\\ $\tau=\sigma\vdash\tau'=\sigma'$ where
$\tau',\sigma'$ are obtained from $\tau,\sigma$ by replacing the
variables simultaneously with arbitrary terms.

We note that the first three axioms are an axiom system for Boolean
algebras, see \cite[Problem 1.1, p.245]{HMTII} (this problem was
solved affirmatively by a theorem prover program).

Consider the four ``logics" (or inference systems) $\vd, \vr$,
[S5,S5,S5], equational logic with (B1 - D5) introduced so far. We
claim that they are equivalent to each other, hence our theorems
stated below apply to all of them.\bigskip

Having formulated our logic $\Ld_3$ in several different ways, we
now formulate some theorems. The first theorem says that this simple
logic $\Ld_3$ is strong enough for doing all of mathematics in it.
It says that we can do set theory in $\Ld_3$ as follows: in place of
formulas $\f$ of set theory we use their ``translated" versions
$\Tr(\f)$ in $\Ld_3$, and then we use the proof system $\vd$ of
$\Ld_3$ between the translated formulas in place of the proof system
of FOL between the original formulas of set theory. Moreover, for
sentences $\f$ in the language of set theory, $\f$ and $\Tr(\f)$
mean the same thing (are equivalent) modulo a ``bridge" $\Delta$
between the two languages. We need this bridge because the language
$\L_{\omega}$ of set theory contains only one binary relation symbol
$\In$ and equality, and the language $\Ld_3$ contains only one
ternary relation symbol $P$. When $f:A\to B$ is a function and
$X\subseteq A$ then $f(X)\de\{ f(a) : a\in X\}$ denotes the image of
$X$ under this function $f$.

\begin{thm}\label{zfc-t}{\sf (Formalizability of set theory in $\Ld_3$)\/}
There is a recursive translation function $\Tr$ from the language
$\L_{\omega}$ of set theory into $\Ld_3$ for which the following are
true for all sentences $\f$ in $\L_{\omega}$:
\begin{itemize}
\item[(i)]
$\ZFC\models\f$\quad iff\quad $\Tr(\ZFC)\vd\Tr(\f)$.
\item[(ii)]
$\ZFC+\Delta\models \f\leftrightarrow\Tr(\f)$,\quad where\\
$\Delta\de\forall xyz[(P(x,y,z)\leftrightarrow(x=y=z\,\,\lor
\,\,\In(x,y))]$.
\end{itemize}
\end{thm}

Theorem~\ref{zfc-t} is proved in section \ref{set-sec}.

\goodbreak
The next theorem is a partial completeness theorem for $\Ld_3$. It
is as good as it can be, see below.

\begin{thm}\label{comp-t}{\sf (Partial completeness theorem for $\Ld_3$)\/}
Let $\L$ be a FOL-language having countably many relation symbols of
each finite arity. There is a recursive subset $K\subseteq\Fmd_3$
and there is a recursive function $\tr$ mapping all $\L$-formulas
into $K$ such that the following are true:
\begin{itemize}
\item[(i)]
$\models\f$\quad iff\quad $\vd\f$ \quad for all $\f\in K$.
\item[(ii)]
$\models\f$\quad iff\quad $\models\tr(\f)$\quad for all
$\L$-sentences $\f$.
\end{itemize}
\end{thm}

According to the above theorem, the proof system $\vd$ is complete
within $K$. But is $K$ big enough? Yes, we can prove any valid
FOL-formula $\f$ by translating it into $K$ and then proving the
translated formula by $\vd$. We know that $\vd$ is not strong enough
to prove all valid $\Fmd_3$ formulas (i.e., $K$ is necessarily a
proper subset of $\Fmd_3$), because as stated in the introduction,
no finite Hilbert-style axiom system can be sound and complete at
the same time for $\Ld_3$. However, we can formulate each sentence
in a slightly different form, namely as $\tr(\f)$ so that this
``version" of $\f$ can now be proved by $\vd$ iff it is valid.

\begin{thm}\label{incomp-t}{\sf (G\"odel style incompleteness theorem for $\Ld_3$)\/}
There is a formula $\f\in\Fmd_3$ such that no consistent recursive
extension $T$ of $\f$ is complete, and moreover, no recursive
extension of $\f$ separates the consequences of $\f$ from the
$\f$-refutable sentences.
\end{thm}

\begin{discussion}\label{discussion}
In Theorems~\ref{zfc-t}-\ref{incomp-t}, at least one at least
ternary relation symbol $P$ is needed in the ``target-language"
$\Ld_3$, the axiom of commutativity ((7)) is needed in the proof
system $\vd$ (because omitting ((7)) from $\vd$ results decidability
of the so obtained proof system, see \cite{NDis}). We do not know
whether complementedness of the closure operators ((6)) is needed or
not. Also, two variables do not suffice because the satisfiability
problem of the two-variable fragment of FOL is decidable.\qed
\end{discussion}

\section{Finding \QRA-reducts in $\Df_3$}\label{qra-sec}

In this section we begin the proof of Theorem~\ref{zfc-t}. For the
definitions of relation algebras, quasi-projective and representable
relation algebras see \cite{TG} or \cite[sec.5.3]{HMTII}. We briefly
recall these. Relation algebras, $\RA$s are Boolean algebras with
operators $\langle A,+,-,;,\inv,1'\rangle$ such that the operators
form an involuted monoid satisfying a further equation. Here,
$\langle A,+,-\rangle$ is the Boolean reduct of the $\RA$ in
question, and $;,\inv,1'$ stand for relation composition, converse,
and identity constant, respectively. The elements $p,q$ in a
relation algebra are called {\it quasi-projections} if
$p\inv;p+q\inv;q\le 1'$ and $p\inv;q=1'+-1'$, and a relation algebra
is called a quasi-projective relation algebra, a $\QRA$, if there is
a pair of quasi-projections in it. We call a relation algebra {\it
representable} if its elements are binary relations and the
operations are union, complementation (w.r.t.\ the biggest element),
relation composition of binary relations, converse of a binary
relation, and the identity relation, respectively (more precisely,
an $\RA$ is representable if it is isomorphic to such a concrete
algebra). Quasi-projective relation algebras are representable, by a
theorem of Tarski.

We show that every $\Df_3$ contains lots of quasi-projective
relation algebras in them. We do this by defining relation algebra
type operations in the term language of $\Df_3$ and proving that
these operations form $\QRA$s in appropriate relativizations. Since
$\QRA$s are representable, this will amount to a ``partial"
representation theorem for $\Df_3$s, and to ``partial" completeness
theorem for $\Ld_3$ (see Thm.\ref{comp-t}), in the spirit of
\cite{HbPhL}. We will work in $\Ld_3$ in place of $\Df_3$.

There will be parameters in the definitions to come.  These will be
formulas in $\Fmd_3$, namely $\delta_{xy}, \delta_{xz}$ with free
variables $\{x,y\}$ and $\{x,z\}$ respectively, together with two
other formulas $p_0, p_1$ with free variables $\{x,y\}$. Thus, if
you choose $\delta_{xy},\delta_{xz},p_0,p_1$ with the above
specified free variables then you will arrive at a \QRA-reduct of
any $\Df_3$ corresponding to these. We get the \QRA-reduct by
assuming some properties of the meanings of these formulas, this
will be expressed by a formula $\Ax$. In section~\ref{set-sec} then
we will choose these parameters so that they fit set theory, which
means that the formula $\Ax$ built up from them is provable in set
theory. Intuitively, the formulas $\delta_{xy},\delta_{xz}$ stand
for equality $x=y,x=z$ and $p_0,p_1$ will be arbitrary pairing
functions.

So, choose formulas $\delta_{xy},\delta_{xz},p_0,p_1$ with the above
specified free variables arbitrarily, they will be parameters of the
definitions to come. To simplify notation, we will {\it not}
indicate these parameters.

We now set ourselves to defining the above relation algebra type
operations on $\Fmd_3$. To help readability, we often write just
comma in place of conjunction in formulas, especially when they
begin with a quantifier. E.g., we write $\exists x(\f,\psi)$ in
place of $\exists x(\f\land\psi)$. Further, $\True$ denotes a
provably true formula, say $\True\de \delta_{xy}\lor\lnot
\delta_{xy}$. First we introduce notation to support the intuitive
meaning of the parameters $\delta_{xy}, \delta_{xz}$ as equality.

\begin{defn}\label{veq-d}(Simulating equality between variables)\\
$\exy \de \delta_{xy}$,\\
$\exz \de \delta_{xz}$,\\
$\eyz \de \exists x (\exy, \exz)$,\\
$\eyx \de \exy$,\\
$\ezx \de \exz$,\\
$\ezy \de \eyz$,\\
$\exx \de \True$,\\
$\eyy \de \True$,\\
$\ezz \de \True$.\qed\\
\end{defn}

\begin{defn}\label{vsubs-d}(Simulating substitution with (simulated) equality)\\
$\varphi\(x,y\) \de \f$,\\
$\varphi\(x,z\) \de \exists y (\eyz , \f)$,\\
$\varphi\(y,z\) \de \exists x (\exy , \f\(x,z\))$,\\
$\varphi\(y,x\) \de \exists z (\exz , \f\(y,z\))$,\\
$\varphi\(z,x\) \de \exists y (\eyz , \f\(y,x\))$,\\
$\varphi\(z,y\) \de \exists x (\exz , \f)$,\\
$\varphi\(x,x\) \de \exists y (\exy , \f)$,\\
$\varphi\(y,y\) \de \exists x (\exy , \f)$,\\
$\varphi\(z,z\) \de \exists x (\exz , \f\(x,x\))$.\qed\\
\end{defn}

\begin{remark}\label{subs-r}
In FOL,  $\f\(u,v\)$ is semantically equivalent with the formula we
get from $\f$ by replacing $x,y$ with $u,v$ everywhere
simultaneously, when $\delta_{xy},\delta_{xz}$ are $x=y, x=z$
respectively.
This is Tarski's fabulous trick to simulate substitutions. 
\qed
\end{remark}

Next we introduce notation supporting intuition about the pairing
functions $p_0, p_1$. First we define some auxiliary formulas. We
will use the notation $2=\{ 0,1\}$, to make the text shorter. Let
$2^*$ denote the set of all finite sequences of $0,1$ including the
empty sequence $\esec$ as well.  If $i,j \in 2^*$ then $ij$ denotes
their ``concatenation'' usually denoted by $i^\cap \!j $, and $|i|$
denotes the ``length'' of~$i$. Further, if $k \in 2$, then we write
$k$ instead of  $\langle k\rangle$ for the sequence $\langle
k\rangle$ of length~$1$. Accordingly, $00$ denotes the sequence
$\langle 0,0\rangle$.

We are going to define $\Fmd_3$-formulas $u_i \eq v_j$ for $u,v\in\{
x,y,z\}$ and  $i, j \in 2^*$. The intuitive meaning of $u_{i_0 \dots
i_n} \eq v_{j_0 \dots j_k}$ is that if $p_0, p_1$ are partial
functions then $p_{i_n} \dots p_{i_0} u = p_{j_k} \dots p_{j_0} v$.
As usual in the partial algebra literature, the equality holds if
both sides are defined and are equal. E.g., the intuitive meaning of
$x_0\eq y_{01}$ is that all of $p_0x, p_0y, p_1p_0y$ exist and
$p_0x=p_1p_0y$.

\begin{defn}\label{pid-d}(Simulating projections)\\
Let $\{ u,v,w\}=\{ x,y,z\}$, $i,j\in 2^*$ and $k\in 2$.\\
$(u_{\esec}\eq v_{\esec}) \de \euv$,\\
$(u_k\eq v_{\esec}) \de p_k\(u,v\)$,\\
$(u_{ik}\eq v_{\esec}) \de \exists w (u_i\eq w_{\esec}, p_k\(w,v\))$\quad if $i\ne \esec$,\\
$(u_i\eq v_j) \de \exists w (u_i\eq w_{\esec},v_j\eq w_{\esec})$\quad if $j\ne \esec$,\\
$(x_i\eq x_j) \de \exists y (\exy, x_i\eq y_j)$,\\
$(y_i\eq y_j) \de \exists x (\exy, x_i\eq y_j)$,\\
$(z_i\eq z_j) \de \exists x (z_i\eq x_{\esec},z_j\eq x_{\esec})$.\qed\\
\end{defn}

We will omit the index $\esec$ in formulas $u_i\eq v_j$, i.e., we
write $u_i \eq v$ and $u \eq v_i$ for $u_i \eq v_{\esec}$ and
$u_{\esec} \eq v_i$ respectively if $i \in 2^*$.

So far we did nothing but introduced notation supporting the
intuitive meanings of the parameters $\delta_{xy}, \delta_{xz}, p_0,
p_1$ as equality and partial pairing functions. Almost any of the
concrete formulas supporting this would do, we only had to fix one
of them since our proof system $\vd$ is very weak, it would not
prove equivalence of most of the semantically equivalent forms. Now
we write up a statement $\Ax$ about the parameters using the just
introduced notation. Let $H\de \{ i\in 2^* : |i|\le 3\}$. Notice
that $H$ is finite.

\begin{defn}[pairing axiom \Ax]\label{ax-d}
We define $\Ax\in\Fmd_3$ to be the conjunction of the union of the
following finite sets (A1),...,(A4) of formulas:
\begin{itemize}
\item[(A1)]
$\{ u_i\eq v_j, v_j\eq w_k \to u_i\eq w_k : u,v,w\in\{x,y,z\},
i,j,k\in H\}$
\item[(A2)]
$\{ u_i\eq v_j, u_{ik}\eq u_{ik}\to u_{ik}\eq v_{jk} :
u,v\in\{x,y,z\}, ik,jk\in H, k\in 2\}$
\item[(A3)]
$\{ u_i\eq u_i, v_j\eq v_j\to \exists w(w_0\eq u_i, w_1\eq v_j) :
u,v,w\in\{x,y,z\}, w\notin\{ u,v\}, i,j\in H\}$
\item[(A4)]
$\{ \exists w\,u\eq w : u,w\in\{ x,y,z\}\}$.\qed
\end{itemize}
\end{defn}

In the above definition,  (A1), (A2), (A4) express usual properties
of the equality, while (A3) states the existence of pairs. We say
that $x$ is a pair if both $p_0$ and $p_1$ are defined on $x$ and
then we think of $x$ as the pair $\langle p_0(x),p_1(x)\rangle$.
That $p_i$ is defined on $x$ is expressed by $p_i(x)\eq p_i(x)$,
i.e., by $x_i\eq x_i$ (for $i\in 2$). Following \cite{TG}, we do not
require pairs to be unique, i.e., for different $u,v$ it can happen
that $u_0=v_0, u_1=v_1$. (This is why $\QRA$s are called {\it
quasi-projective} $\RA$s and not just projective $\RA$s in
\cite{TG}.) In the next section, just for simplicity, we will use a
stronger axiom $\SAx$ in place of $\Ax$ in which we require
uniqueness of pairs.

We are ready to define our relation-algebra type operations on
$\Fmd_3$. They will have the intended meanings on formulas with one
free variable $x$, where  $x$ denotes a pair.  This is expressed by
the definition of $\dra$, the universe of the algebra defined below.
If we assume uniqueness of pairs (as in $\SAx$ later) then the
definition of $\dra$ in Def.\ref{qra-d} below can be simplified to
be $\dra\de\{\f\in\Fmd_3^1 :\Ax\vd\f\to\pair\}$, where $\pair$ is
the formula expressing that $x$ is a pair. Since $\f$ has only one
free variable $x$ which is a pair, we can think of $\f$ as a unary
relation of pairs, i.e., as a binary relation. With this intuition,
the definitions of the operations $\comp,\conv,\id$ below in
Def.\ref{qra-d} are the natural ones, see Figure~1. For more on the
intuition behind Def.\ref{qra-d} see the remark after the
definition.

\begin{defn}[relation algebra reduct $\Dra$ of $\Fmd_3$]
\label{qra-d} Let $\f,\psi\in\Fmd_3$.
\begin{itemize}
\item[]
$\pair \de \exists y p_0 \land \exists y p_1$. \\
$\f u_i \de \exists x(x\eq u_i,\f)$ \quad if $u\in\{ y,z\}$ and
$i\in 2^*$.
\item[]
$\f\comp\psi \de \exists y(\f y_0, \psi y_1, x_0\eq y_{00},
y_{01}\eq y_{10}, y_{11}\eq x_1)$,\quad see Figure~1,\\
$\f\conv \de \exists y(\f y, y_0\eq x_1, y_1\eq x_0)$,\\
$\id \de x_0\eq x_1$,\\
$\div\, \f \de \pair\land\lnot\f$,\\
$\f+\,\psi \de \f\lor\psi$.
\item[]
$\dra \de \{\f\in\Fmd_3 : \Ax\vd\f\leftrightarrow \psi\comp\id\mbox{
 for some  }\psi\in\Fmd_3^1\}$,\\
 $\Dra \de \langle\dra,+,\,\div,\comp,\conv,\id\rangle$.\qed
\end{itemize}
\end{defn}

\begin{figure}[!h]\label{comp-fig}
\centerline{\psfig{file=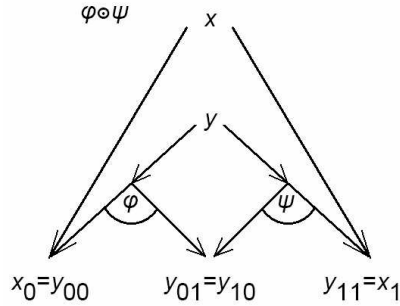,height=45mm}}
\caption{Illustration of $\varphi \comp \psi$}
\end{figure}

Let us define $x\sim y\de x_0\eq y_0, x_1\eq y_1$, and let us call
the pairs $x,y$ equivalent if $x\sim y$. Since we do not require
pairs to be unique in $\Ax$, we may have distinct $x,y$ which are
equivalent. However, in the above definition, we can see that the
result of an operation from $\comp,\conv,\id$ is always closed under
the equivalence relation $\sim$, because, intuitively, the result
depends only on $x_0, x_1$, thus if, say, $\f\comp\psi$ holds at $x$
and $x\sim y$, then $\f\comp\psi$ holds at $y$, too. (Formally,
$\Ax\vd\f\comp\psi\land x\sim y\to (\f\comp\psi)y$.) From this one
can see that $\f\comp\id$ represents the same binary relation as
$\f$ composed with $\sim$; and thus $\dra$ consists of those
formulas which do not distinguish equivalent pairs.

With the intuition that $\f$ represents a binary relation (coded as
a unary relation on pairs), we could have defined, say,
$\f\comp\psi$ as $\exists yz(\f y, \psi z, x_0\eq y_0, y_1\eq z_0,
z_1\eq x_1)$. This definition would leave us with one register (or
free variable) to work with, namely $x$ (because $x_0, x_1$ are
recoverable from $y,z$). It is more convenient to code up every
relevant data in one register ($y$ in Def.\ref{qra-d}) and so have
two registers (namely, $x,z$) to work with. This is how we defined
the relation algebraic operations in Def.\ref{qra-d} above.

\goodbreak
The next theorem is the heart of formalizing set theory in $\Fmd_3$.
Let us define the equivalence relation $\equiv_{T}$ on $\Fmd_3$ by
$$ \f\equiv_{T}\psi \de T\vd\f\leftrightarrow\psi,$$ where
$T\subseteq\Fmd_3$. Note that with using this notation we have
$$ \dra = \bigcup\{ \f\comp\id\slash\equiv_{\Ax}\ \  :\
\f\in\Fmd^1_3\}.$$

\goodbreak
\begin{thm}\label{dra-t}
\begin{itemize}
\item[(i)]
$\Dra$ is an algebra, i.e., the operations $+,\div,\comp,\conv,\id$
do not lead out of $\dra$.
\item[(ii)]
$\dra\supseteq\{ \f\comp\psi : \f,\psi\in\Fmd_3^1\}$.
\item[(iii)]
$\Dra\slash\equiv_{\Ax}$ is a relation algebra.
\item[(iv)]
The images of the formulas $x_1\eq x_{00}$ and $x_1\eq x_{01}$ form
a quasi-projection pair in $\Dra\slash\equiv_{\Ax}$.
\end{itemize}
\end{thm}

\noindent{\bf Proof.} The proof of the analogous theorem in
\cite[Thm.9, p.43]{NDis} goes through with some modifications. We
indicate here these modifications.

The ``proof explanations" (CA),\dots, (KV) introduced in
\cite[p.44]{NDis} can be used in our proof, too, except that we
always have to check whether the explanation uses only the Df-part
of (CA).  If it uses axiom C7 of $\CA$, then we have to give an
alternate proof. Let (Df) denote the part of (CA) which does not use
C7. We can use (UV) because it can be derived from $\Ax$ and (Df).
Of course, throughout we have to change $=$ to $\eq\,$.

We now go through the proof given in \cite[pp.46-64]{NDis} and
indicate the changes needed for our proof. All the statements on
pp.48-54 beginning with (A1) and ending with (S13) follow from $\Ax$
and (Df). In fact, we could just add these statements to $\Ax$ since
they do not contain formula-schemes denoting arbitrary formulas
(such as, e.g., (0) on p.54 does), and so they amount to finitely
many formulas only. Because of this, we do not indicate the changes
needed in the proofs of these items.

We have to avoid statement (0) at the end of p.54 by all means,
because it uses axiom C7 of (CA) essentially. Fortunately, we do not
really use (0) in the proof, changing $\f$ to $\f x$ in some steps
will suffice for eliminating (0).

The proof of (2\slash a) on p.55 has to be modified, and that can be
done as follows. Recall from \cite{NDis} that $\f$ has at most one
free variable $x$, and $\f u\de \exists x(x\eq u, \f)$ when $u$ is
different from $x$ while $\f x\de\exists y(y\eq x, \f y)$.
\smallskip

\noindent Assume $x\notin\{ u,v\}$. The proof for $\f u, u\eq v\vd
\f v$ is as on \cite[p.55]{NDis}. Then

\bigskip
\noindent $\f x, x\eq z \vd\qquad\mbox{by definition of $\f x$}\\
\exists y(y\eq x, \f y), x\eq z \vd\\
\exists y(y\eq x, \f y, x\eq z )\vd\qquad\mbox{ by \Ax}\\
\exists y(y\eq z, \f y) \vd\qquad\mbox{ by the case when $x\notin\{ u,v\}$ }\\
\exists y(\f z) \vd\\
\f z.$

\bigskip
\noindent $\f x, x\eq y \vd\qquad\mbox{ by \Ax}\\
\exists z(z\eq x, \f x, x\eq y) \vd\qquad\mbox{ by \Ax}\\
\exists z(z\eq y, \f x, x\eq z) \vd\qquad\mbox{ by the previous case}\\
\exists z(z\eq y, \f z)\vd\qquad\mbox{ by the case when $x\notin\{ u,v\}$ }\\
\f y.$
\bigskip

\noindent $\f y, y\eq x \vd\\
\exists y(y\eq x, \f y) \vd\quad\mbox{ by definition of $\f x$}\\
\f x .$
\bigskip

\noindent $\f z, z\eq x \vd\quad\mbox{ by $\Ax$}\\
\exists y(y\eq z, z\eq x, \f z) \vd\quad\mbox{ by $\Ax$}\\
\exists y(y\eq x, z\eq y, \f z) \vd \quad\mbox{ by the case when $x\notin\{ u,v\}$ }\\
\exists y(y\eq x, \f y) \vd\\
\f x.$
\bigskip

On p.58, in the last line of the proof of (9) we have to write
$\gamma x$ in place of $\gamma$. Similarly, on p.59, in lines 7 and
8 we have to change $\psi$ to $\psi x$ and then we can cross out
reference to (0).

In the last line of the proof of (18) on p.62 we have to use the
following, which is practically C7 for formulas of form
$\f\comp\psi$:

\begin{itemize}
\item[(C)]
$[\lnot(\f\comp\psi)]z\leftrightarrow \lnot[(\f\comp\psi)z]\ .$
\end{itemize}

\noindent Recall that $\Delta(x,y)$ denotes the following formula:
$x_0\eq y_{00}, y_{01}\eq y_{10}, x_1\eq y_{11}$ (cf.,
Figure~\ref{comp-fig}). To prove (C) we will prove the following two
statements from which it follows immediately.

\begin{itemize}
\item[(C1)]
$[\lnot(\f\comp\psi)]z\leftrightarrow \neg\exists y(\Delta(z,y),\f
y_0,\psi y_1)$ .
\item[(C2)]
$(\f\comp\psi)z\leftrightarrow \exists y(\Delta(z,y),\f y_0,\psi
y_1)$ .
\end{itemize}

\noindent {\it Proof of (C2):}\\
\noindent $(\f\comp\psi)z\lr\quad\mbox{ by definition of $\chi z$}\\
\exists x(x\eq z, \f\comp\psi) \lr\quad\mbox{ by definition of $\comp$}\\
\exists x(x\eq z, \exists y(\Delta(x,y),\f y_0, \psi y_1))\lr\quad\mbox{ by SZV}\\
\exists x\exists y(x\eq z,\Delta(x,y),\f y_0, \psi y_1)\lr\quad\mbox{ by $\Ax$}\\
\exists x\exists y(x\eq z,\Delta(z,y),\f y_0, \psi y_1)\lr\quad\mbox{ by SZV}\\
\exists x(x\eq z,\exists y(\Delta(z,y),\f y_0, \psi y_1))\lr\quad\mbox{ by SZV}\\
\exists x(x\eq z),\exists y(\Delta(z,y),\f y_0, \psi y_1)\lr\quad\mbox{ by $\Ax$}\\
\exists y(\Delta(z,y),\f y_0, \psi y_1)\quad\mbox{ and we are
done.}$
\bigskip

\noindent {\it Proof of (C1):}\\
\noindent $[\neg(\f\comp\psi)]z\lr\quad\mbox{ by definitions of $\chi z$ and $\comp$}\\
\exists x(x\eq z, \neg\exists y(\Delta(x,y),\f y_0, \psi y_1))\lr\quad\mbox{ by SZV}\\
\exists x(\forall y(x\eq z),\forall y\neg(\Delta(x,y),\f y_0, \psi y_1))\lr\quad\mbox{ by Df}\\
\exists x\forall y(x\eq z, \neg(\Delta(x,y),\f y_0, \psi y_1))\lr\quad\mbox{ by $\Ax$ and BA}\\
\exists x\forall y(x\eq z, \neg(\Delta(z,y),\f y_0, \psi y_1))\lr\quad\mbox{ by Df, SZV}\\
\exists x(x\eq z,\neg\exists y(\Delta(z,y),\f y_0, \psi y_1))\lr\quad\mbox{ by SZV}\\
\exists x(x\eq z),\neg\exists y(\Delta(z,y),\f y_0, \psi y_1)\lr\quad\mbox{ by $\Ax$}\\
\neg\exists y(\Delta(z,y),\f y_0, \psi y_1)\quad\mbox{ and we are
done}.$
\bigskip

That's all the changes we have to do in the proof given in
\cite[pp.46-64]{NDis}! \qquad {\bf QED}\bigskip

\section{Finding \CA-reducts in $\Df_3$}\label{ca-sec}

Simon~\cite{Simon} defines a $\CA_n$-reduct  in every $\QRA$, for
all $n\in\omega$, and also proves that these reducts are
representable. We will use, in this paper, the $\CA_3$-reduct of our
$\QRA$ defined in the previous section, i.e., we will use the
$\CA_3$-reduct of $\Dra\slash\equiv_{\Ax}$.

We will use the following stronger form of $\Ax$, just for
convenience:

\begin{defn}\label{sax-d}(strong pairing axiom \SAx.)
We define $\SAx\in\Fmd_3$ to be the conjunction of the union of the
finite sets (A1),...,(A4) of Def.\ref{ax-d} together with the
following:
\begin{itemize}
\item[(A5)]
$\{\ x_0\eq y_0, x_1\eq y_1\to x\eq y, \ \ x_0\eq x_0\leftrightarrow
x_1\eq x_1\}$.\qed
\end{itemize}
\end{defn}

The formulas in (A5) express that pairs are unique, and that $p_0$
is defined exactly when $p_1$ is defined. We use (A5) for
convenience only, this way formulas will be shorter. We could omit
(A5) on the expense that formulas will be longer and more
complicated. If we assume $\SAx$ then $\f\land\pair\in\dra$ for all
$\f\in\Fmd_3^1$.

Every $\QRA$ has a $\CA_3$-reduct, which is representable, see Simon
\cite{Simon}. The following definition is recalling this
$\CA_3$-reduct from \cite{Simon} in our special case of
$\Dra\slash\equiv_{\SAx}$. The definition below is simpler than in
\cite{Simon} because we will assume uniqueness of pairs in $\SAx$,
which is not assumed in \cite{Simon}. We will use the abbreviation
\begin{itemize}
\item[]
$\f x_i\de\exists y(y\eq x_i,\f y)$, when $\f\in\Fmd^1_3$ and $i\in
2^*$.
\end{itemize}
Beware: $\vd \f\leftrightarrow\f x$ usually is not true for
$\f\in\Fmd_3^1$. (For the definition of $\f y_j$ see
Def.\ref{qra-d}.)

The intuitive meaning of Def.\ref{ca-d} below is similar to the one
of Def.\ref{qra-d}. The universe of our $\CA_3$ will consist of
those formulas which depend only on $x_1$ such that $x_1$ is a
triplet; we will look at such a formula as representing a set of
these triplets, i.e., a ternary relation. With this is mind, then in
the definition below, $\d_{ij}$ represents the set of those triplets
whose $i$-th and $j$-th components equal, $\T_i$ is the binary
relation on triplets which correlates two triplets iff only their
$i$-th components may differ.

\begin{defn}\label{ca-d}(cylindric reduct $\Dca$ of $\Fmd_3$.)\\
$\Triplet\de x_{11}\eq x_{11}$,\quad see Figure~2,\\
$(0)\de 0,\quad (1)\de 10,\quad (2)\de 11$,\quad and for all $i,j<3$ and $\f,\psi\in\Fmd_3^1$\\
$\T_i\de\Triplet x_0\land\Triplet x_1\land\bigwedge\{ x_{0(j)}\eq
x_{1(j)} : i\ne j<3\}$, \\
$\c_i\f\de\f\comp\T_i$,\\
$\D_{ij}\de \Triplet x_1\land x_{1(i)}\eq x_{1(j)}$,\\
$\mi\, \f \de \Triplet x_1\land\lnot\f$,\\
$\f+\,\psi \de \f\lor\psi$,\\
$\dca\de\{ \f\in\Fmd_3 : \SAx\vd \f\leftrightarrow
(\psi x_1\land\Triplet x_1)\ \  \mbox{ for some } \psi\in\Fmd_3^1\}$,\\
 $\Dca \de \langle\dca,+,\,-,\,\c_i,\D_{ij} : i,j<3\rangle$.\qed
\end{defn}

\begin{figure}\label{triplet-fig}
\begin{center}
\psfrag*{X}[b][b]{$x$} \psfrag*{X0}[t][t]{$x_{(0)}$}
\psfrag*{X1}[t][t]{$x_{(1)}$} \psfrag*{X2}[t][t]{$x_{(2)}$}
\includegraphics[keepaspectratio, width=0.3\textwidth]{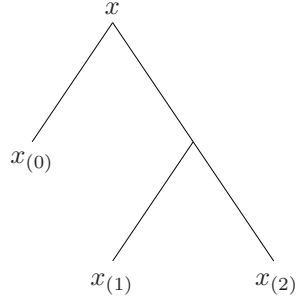}
\caption{$x$ represents the triplet $\langle
x_{(0)},x_{(1)},x_{(2)}\rangle$  when $x_{11}$ is defined}
\end{center}
\end{figure}

\goodbreak
\begin{thm}\label{cared-t} The set $\dca$ is closed under the
operations defined in Def.\ref{ca-d} and the algebra
$\Dca\slash\equiv_{\SAx}$ is a representable $\CA_3$.
\end{thm}

\noindent{\bf Proof.} We show that the algebra
$\Dca\slash\equiv_{\SAx}$ is the $\CA_3$-reduct of the
quasi-projective relation algebra $\Dra\slash\equiv_{\SAx}$, as
defined in \cite{Simon}. Let $\p\de x_1\eq x_{00}$, $\q\de x_1\eq
x_{01}$. In the following we will omit referring to $\equiv_{\SAx}$,
so we will look at $\p$ as an element of $\Dra\slash\equiv_{\SAx}$,
while only $\p\slash\equiv_{\SAx}$ is such. Recall from
Thm.\ref{dra-t}(iv) that $\p,\q$ form a pair of projections in
$\Dca\slash\equiv_{\SAx}$. See Figure~3.

\begin{figure}[!h]
\begin{center}
\psfrag*{0}[r][r]{$0$} \psfrag*{1}[l][l]{$1$}
\psfrag*{0!}[b][b]{$0$} \psfrag*{x}[b][b]{$x$}
\psfrag*{1!}[b][b]{$1$}
\includegraphics[keepaspectratio, width=0.5\textwidth]{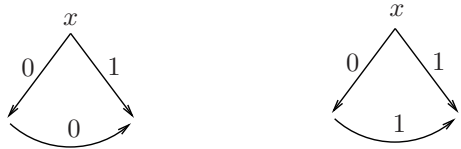}
\caption{$\p$ and $\q$ represent the pairing functions $p_0, p_1$
coded as unary relations of pairs.}
\end{center}
\end{figure}

\goodbreak
Let
\begin{itemize}
\item[]
$\e\de \Triplet x_0\land x_0\eq x_1$,\\
$\pi_{(i)}\de\Triplet x_0\land x_1\eq x_{0(i)}$,\\
$\chi_{(i)}\de\Triplet x_0\land\Triplet x_1\land x_{0(i)}\eq
x_{1(i)}$.
\end{itemize}
\noindent Now, one can show that $\e$ is the same as $\in\!\!^{(3)}$
in \cite[Def.3.1]{Simon}, i.e.,
\begin{itemize}
\item[]
$\SAx\vd\e\leftrightarrow(\q\comp\q\comp\q\conv\comp\q\conv\land\id).$
\end{itemize}
Similarly, one can show that $\pi_{(i)}, \chi_{(i)}$ are the same as
the ones in \cite{Simon}, and $\Triplet$ is the same as $1^{(3)}$ in
\cite[Def.3.1]{Simon}, e.g.,
\begin{itemize}
\item[]
$\SAx\vd\pi_{(1)}\leftrightarrow(\e\comp\q\comp\p)$,\\
$\SAx\vd\chi_{(1)}\leftrightarrow(\pi_{(1)}\comp\pi_{(1)}\conv)$,\\
$\SAx\vd\Triplet\leftrightarrow(\pair\comp\e)$.
\end{itemize}
From this one can show, similarly to the above, that \cite{Simon}'s
$\t_i, \d_{ij},\t$ are the same as our $T_i, \d_{ij}$ and $\e$ in
Def.\ref{ca-d} and above. Thus, our reduct $\Dca$ is the same as the
3-reduct defined in \cite{Simon}, and then we can use
\cite[Thm.3.2,Thm.5.2]{Simon}.\qquad {\bf QED}\bigskip

\section{Formalizing set theory in $\Ld_3$}\label{set-sec}

The 3-variable restricted fragment $\L_3$ of $\L_\omega$ is defined
as follows. Language: three variable symbols, $x,y,z$, one binary
relational symbol $\In$, so there is one atomic formula, namely
$\In(x,y)$. Logical connectives: $\lor, \neg, \exists x, \exists y,
\exists z$ and $u=v$ for $u,v\in\{ x,y,z\}$. We denote the set of
formulas by $\Fm_3$. Derived connectives are $\forall, \lor, \to,
\leftrightarrow$. The proof system $\vdr$ of $\L_3$ is a Hilbert
style one defined in \cite[Part II, p.157]{HMTII}. The word-algebra
of $\L_3$ is denoted as $\gFm_3$, it is the absolutely free
$\CA_3$-type algebra generated by the formula $\In(x,y)$.
$\Fm_\omega$ denotes the set of formulas of $\L_\omega$, and if
$\Fm$ is one of $\Fmd_3, \Fm_3, \Fm_\omega$, then $\Fm^0, \Fm^1,
\Fm^2$ denote the set of formulas in $\Fm$ with no free variables,
with one free variable $x$, and with two free variables $x,y$,
respectively.

In this section we prove Theorem \ref{zfc-t} stated in section
\ref{fmd-sec}. As a first step, we define a translation function
$\h$ from $\Fm_3$ to $\Fmd_3$. This will be analogous to the one
defined in \cite{NDis}, but here a novelty is that the vocabulary of
$\Fm_3$ is different (disjoint) from that of $\Fmd_3$, and we will
have to pay attention to this difference. Namely, $\Fm_3$ contains
one binary relation symbol $\In$ and equality $u=v$ for $u,v\in\{
x,y,z\}$ while $\Fmd_3$ contains one ternary relation symbol $P$ and
does not contain equality. The formula $\Delta$ introduced in
section \ref{fmd-sec} bridges the difference between the two
languages. Namely, $\Delta$ is stated in a language which contains
both $\Fm_3$ and $\Fmd_3$ and $\Delta$ is a definition of $P$ in
$\Fm_3$, so it leads from $\Fm_3$ to $\Fmd_3$. But $\Delta$ is a
two-way bridge, because of the following. Let
\begin{align*}
\Delta'\de [\In(x,y)\leftrightarrow \forall
zP(x,y,z)]\land[x=y=z\leftrightarrow (P(x,y,z)\land\lnot\forall
zP(x,y,z))]\ .
\end{align*}
Then $\Delta'$ provides a definition of equality $=$ and $\In$ in
$\Fmd_3$, and thus it is a bridge leading from $\Fmd_3$ to $\Fm_3$.
Now, the two definitions are equivalent in non-trivial models,
namely
$$ \exists xy(x\ne y)\models \Delta\leftrightarrow \Delta'\ .$$

As usual, we begin with definitions. To start, we work in $\Fm_3$.
Recall the concrete definitions of $p_0, p_1$, and $\pi$ from
\cite[p.71, p.35]{NDis}. Define

\bigskip
$\pi^+\de \pi\land\\
\phantom{\pi^+\de\quad\ \  }\forall xy[\exists z(\p_0\nz
x,z\bz,\p_0\nz y,z\bz),\exists
z(\p_1\nz x,z\bz,\p_1\nz y,z\bz)\to x=y]\land\\
\phantom{\pi^+\de\quad\ \  }\forall x(\exists
y\p_0(x,y)\leftrightarrow \exists y\p_1(x,y))\land\exists xy(x\ne
y).$
\bigskip

\noindent Then $p_0, p_1, \pi, \pi^+$ are formulas of $\Fm_3$. We
note that the notation $\f\nz u,v\bz$ in \cite{NDis} is the same as
our $\f\(u,v\)$ introduced in Def.\ref{vsubs-d}, except that instead
of $x\eq y$ etc.\ the notation $\f\nz u,v\bz$ uses real equality
$x=y$ etc.\ available in $\Fm_3$.

Next, we work in $\Fmd_3$ and we get the versions of these formulas
that the definition $\Delta$ of $P$ provides. We will write out the
details. First we fix the parameters $\delta_{xy}, \delta_{xz}, p_0,
p_1$ occurring in the formula $\SAx$ of $\Fmd_3$.

\begin{defn}\label{E-d}(fixing the parameters $\delta_{xy}, \delta_{xz}$ of $\Fmd_3$)\\
$E\de \forall z\,P(x,y,z)$,\\
$D\de P(x,y,z)\land\lnot E$,\\
$\delta_{xy}\de \exists zD$,\\
$\delta_{xz}\de \exists yD$.\qed
\end{defn}

The next definition fixes the parameters $p_0,p_1$ of $\Fmd_3$. To
distinguish them from their $\Fm_3$-versions, we will denote them by
$\p_0, \p_1$. The definition below is a repetition of the definition
of the corresponding formulas $p_0,p_1$ on p.71 of \cite{NDis} such
that we write $E$ and the above defined concrete $\delta_{xy},
\delta_{xz}$ in place of $\In$ and $x=y$, $x=z$. We will use the
notation introduced in section \ref{fmd-sec} and we will use
notation to support set theoretic intuition. Thus, if we introduce a
formula $\f$ denoted as, say, $x\eq\{ y\}$, then $u\eq\{ v\}$
denotes the formula $\f\(u,v\)$ (see
Definitions~\ref{vsubs-d},\ref{veq-d}). Below, ``op" abbreviates
``ordered pair" (to distinguish it from the formula $\pair$ defined
earlier.).

\begin{defn}\label{p-d}(fixing the parameters $\p_0,\p_1$ of $\Fmd_3$)\\
$u\inn v\de E\(u,v\)$,\quad for $u,v\in\{ x,y,z\}$,\\
$x\eq\{ y\}\de \forall z(z\inn x \leftrightarrow z\eq y)$,\\
$\{ x\}\inn y\de \exists z(z\eq\{ x\}, z\inn y)$,\\
$x\eq\{\{ y\}\}\de \exists z(z\eq\{ y\}, x\eq\{ z\})$,\\
$x\inn\cup y\de \exists z(x\inn z, z\inn y)$,\\
$\op(x)\de \exists y\forall z(\{ z\}\inn x\leftrightarrow y\eq
z)\land$\\
$\phantom{\op(x)\de} \forall yz[(y\inn\cup x, \lnot\{ y\}\inn x,
z\inn\cup x, \lnot\{ z\}\inn x)\rightarrow y\eq z], \forall y\exists
z(y\inn x\rightarrow z\inn y)$,\\
$\p_0\de \op(x)\land \{ y\}\inn x$,\\
$\p_1\de \op(x)\land[x\eq\{\{ y\}\}\lor (y\inn\cup x, \lnot\{
y\}\inn x)]$.\qed\\
\end{defn}

It is not hard to check that
\begin{align}\label{bridge}
\Delta,\exists xy(x\ne y)\models\ \pi\!\leftrightarrow\Ax,\
\pi^+\!\!\leftrightarrow\SAx .
\end{align}

We are ready to define our translation mapping $\h$  from $\Fm_3$ to
$\Fmd_3$. For a formula $\f\in\Fmd_3^2$ and $i,j\in 2^*$ we define
$$\f\(x_i,x_j\)\de\exists yz(y\eq x_i, z\eq x_j, \f\(y,z\)).$$

\begin{defn}\label{h-d}(translation mapping $\h$)
\begin{itemize}
\item[(i)]
$\h': \Fm_3\to \Fmd_3^1$ is defined by the following:\\
$\h'$ is a homomorphism from $\gFm_3$ into $\Dca$ such that\\
$\h'(\In(x,y))\de E\(x_{1(0)},x_{1(1)}\)$.
\item[(ii)]
$\SAx^*\de\SAx\land\forall x(\Triplet x_1\to\h'(\pi^+))$.
\item[(iii)]
We define the mapping $\h:\Fm_3\to\Fmd_3^0$ as\\
$$\h(\f)\de \forall x([\SAx^*\land\Triplet x_1]\to\h'(\f)).\qed$$
\end{itemize}
\end{defn}

We say that a translation function $f$ is {\it Boolean preserving
w.r.t.\ $\vd$} iff for all sentences $\f,\psi\in\Fmd_3$ we have that
$\vd f(\f\land\psi)\leftrightarrow(f(\f)\land f(\psi))$ and $\vd
f(\f\to\psi)\to(f(\f)\to f(\psi)).$

\begin{thm}\label{h-t} Let $\f$ be a sentence
and $T$ be a set of sentences of $\Fm_3$. Then the following
(i)-(iii) are true.
\begin{itemize}
\item[(i)]
$T\cup\{\pi^+\}\models\f\ \Leftrightarrow\ \h(T)\vd\h\f$.
\item[(ii)]
$\pi^+\land\Delta\models \f\leftrightarrow\h(\f)$.
\item[(iii)]
$\h$ is Boolean preserving and $\SAx^*\vd\h(\lnot\f)\to\lnot\h(\f)$.
\end{itemize}
\end{thm}

\noindent{\bf Proof.} \underbar{Proof of (ii):}  Let $\Mm=\langle
M,\In^{\Mm}, P^{\Mm}\rangle$ be a model of $\pi^+\land\Delta$. Let
$V\de\{ x,y,z\}$ and $\Val\de {}^VM$, the set of evaluations of the
variables in $\Mm$. Now, by $\Mm\models\pi^+\land\Delta$ we have
that $p_i$ and $\p_i$ have the same meanings in $\Mm$, and they form
pairing functions. Then one can show by induction that
\begin{align}
\Mm\models (u_i=v_j)[k]\quad\mbox{iff}\quad\Mm\models (u_i\eq
v_j)[k]\quad\mbox{iff}\quad k(u)_i=k(v)_j \mbox{ in }\Mm.
\end{align}
We say that $a\in M$ is a {\it triplet} iff $a_{11}$ is defined.
Then $a_0, a_{10}$ are also defined by $\Mm\models\pi^+$. If $a$ is
a triplet, then we assign an evaluation $\val(a)\in\Val$ to $a$ such
that $\val(a)$ assigns to $x,y,z$ the elements $a_{(0)}, a_{(1)},
a_{(2)}$ respectively. One can prove by induction the following
statement: For all $\f\in\Fm_3$ and $k\in\Val$ we have
\begin{align}\label{hval}
\Mm\models\h'(\f)[k]\quad\mbox{iff}\quad \big(k(x)_1\mbox{ is a
triplet and }\Mm\models\f[\val(k(x)_1)]\big).
\end{align}
From (\ref{hval}) we get the following:
\begin{align}\label{hcon}
\Mm\models \f\leftrightarrow \forall x(\Triplet
x_1\to\h'\f)\quad\mbox{for all }\quad \f\in\Fm_3^0.
\end{align}
Indeed, assume $\Mm\models\f$, then $\Mm\models\f[k]$ for all $k\in
\Val$. We show $\Mm\models\forall x(\Triplet x_1\to\h'\f)$. Indeed,
let $k\in\Val$ be such that $k(x)_1$ is a triplet. By $\Mm\models\f$
then $\Mm\models\f[\val(k(x)_1)]$, by (\ref{hval}) then
$\Mm\models\h'\f[k]$ and we are done.
Assume next $\Mm\notmodels\f$, then $\Mm\notmodels\f[k]$ for all
$k\in\Val$ because $\f$ is a sentence (i.e., does not contain free
variables). We show $\Mm\notmodels\forall x(\Triplet x_1\to\h'\f)$.
Indeed, let $a\in M$ be such that $a_{111}$ is defined. Such an
$a\in M$ exists by $\Mm\models\pi$. Let $k\in\Val$ be such that
$k(x)=a$. Then $k(x)_1$ is a triplet and
$\Mm\notmodels\f[\val(k(x)_1]$, so $\Mm\notmodels\h'\f[k]$ by
(\ref{hval}), so $\Mm\notmodels(\Triplet x_1\to\h'\f)[k]$, so
$\Mm\notmodels\forall x(\Triplet x_1\to\h'\f)$. This shows that
(\ref{hcon}) indeed holds.

From (\ref{hcon}) and $\Mm\models\pi^+\land\Delta$ then we have
$\Mm\models\SAx^*$. This together with (\ref{hcon}) and the
definition of $\h$  finishes the proof of (ii).

\bigskip
\noindent\underbar{Proof of (iii):} The proofs of (4) and (6) in
\cite[p.73]{NDis}, which prove that $\kappa$ is Boolean-preserving
w.r.t.\ $\vdr$, work for showing that $\h$ is Boolean preserving
w.r.t.\ $\vd$, because $\h$ has the same ``structure" as $\kappa$.
Similarly, the proof of (5) in \cite[p.73]{NDis} is good for proving
the second statement of the present (iii).\bigskip

\noindent\underbar{Proof of (i):} First we prove (i) for the special
case when $T$ is the empty set. Let $\f$ be a sentence of $\Fm_3$,
we want to prove $\notvd\h(\f)$ implies $\pi^+\notmodels\f$. So,
assume that $\notvd\h\f$,\ \  i.e., $\notvd\forall
x(\SAx^*\land\Triplet x_1\to\h'\f)$. Then\ \ $\SAx\notvd\forall
x(\Triplet x_1\land\h'(\pi^+)\to\h'(\f))$, \ \ i.e.,\ \
\begin{align}\label{notvd}
\SAx\notvd\h'(\pi^+\!\!\to\f). \end{align} Recall that $\h'$ is a
homomorphism from $\gFm_3$ to $\Dca\slash\equiv_{\SAx}$, and the
latter is a representable $\CA_3$. Let $\psi\de\pi^+\!\!\to\f$. By
(\ref{notvd}) we have that the image of $\psi$ is not $1$ under
$\h'$, therefore there is a homomorphism $g$ from
$\Dca\slash\equiv_{\SAx}$ to a cylindric set algebra $\Cc$ such that
the image of $\h'\psi$ under $g$ is not $1$. Let $f\de g\h'$, then
\begin{align}\label{notmodel}
f:\gFm_3\to\Cc\qquad \mbox{ and }\qquad  f(\psi)\neq 1.
\end{align}
Let $U$ be the base set of $\Cc$, let $R\de \{\langle s_0,
s_1\rangle : s\in f(\In(x,y))\}$ and define the model $\Mm$ as
$\langle U,R\rangle$. Then for all $\f\in\Fm_3$ and $s\in U^3$ we
have that $\Mm\models\f[s]$ iff $s\in f(\f)$. Thus
$\Mm\notmodels\psi$ by (\ref{notmodel}), and we are done with
showing $\pi^+\notmodels\f$.

In the other direction, we have to show that $\vd\h\f$ implies
$\pi^+\models\f$. By soundness of the proof system $\vd$ we have
$\models\h\f$, then by (ii) we have $\pi^+\land\Delta\models\f$.
Since $P$ does not occur in $\f$ and in $\pi^+$, this means that
$\pi^+\models\f$ and we are done.

Next, assume that $T$ is a set of sentences of $\Fm_3$. We want to
show that $T\cup\{\pi^+\}\models\f$ iff $\h(T)\vd\h(\f)$. Now,
\bigskip

\noindent
$T\cup\{\pi^+\}\models\f$ iff (by compactness of $\Fm_3$)\\
$T_0\cup\{\pi^+\}\models\f$ for some finite $T_0\subseteq T$,  iff\\
$\pi^+\models\bigwedge T_0\to\f$ for some finite $T_0\subseteq T$,  iff(by first part of (i))\\
$\vd\h(\bigwedge T_0\to\f)$ for some finite $T_0\subseteq T$. Then,
by Boolean preserving of $\h$ \\
$\vd\bigwedge\h(T_0)\to\h(\f)$, and so by Modus Ponens we get\\
$h(T_0)\vd\h(\f)$. Conversely, from this we get by the soundness of
$\vd$ that\\
\begin{align}\label{mod}
h(T_0)\models\h(\f). \end{align}
 From here on we have to deal with
the difference between the two languages.

Let $\Mm$ be an arbitrary model of $\L_3$, so $\Mm$ contains one
binary relation, say $\Mm=\langle M,\In^{\Mm}\rangle$. Define
$P^{\Mm}$ according to the definition $\Delta$, i.e., $P^{\Mm}\de\{
\langle a,a,a\rangle : a\in M\}\cup\{\langle a,b,c\rangle\in M\times
M\times M : \langle a,b\rangle\in\In^{\Mm}\}$. Let $\Mm^+\de\langle
M,\In^{\Mm},P^{\Mm}\rangle$ be the expansion of $\Mm$ with this new
relation, and let $\Mm^-\de\langle M,P^{\Mm}\rangle$ be the reduct
of this expansion to the language $\Ld_3$.

Assume now $\Mm\models T_0\cup\{\pi^+\}$. Then \\
$\Mm^+\models T_0\cup\{\pi^+\land\Delta\}$. Thus by (ii) we have
that\\
$\Mm^+\models\h(T_0)$, then \\
$\Mm^-\models\h(\f)$ by (\ref{mod}). Thus\\
$\Mm^+\models \h\f\land\pi^+\land\Delta$, so by (ii)\\
$\Mm^+\models \f$,  and so\\
$\Mm\models\f$  and thus, since $\Mm$ was an arbitrary model of $\L_3$\\
$\T_0\cup\{\pi^+\}\models\f$ and we are done.\qquad{\bf QED}\bigskip

We are almost done with proving Thm.\ref{zfc-t}, all we have to do
is to use a connection between $\L_\omega$ and $\L_3$ established in
\cite{NDis}, \cite{NPrep}. \bigskip

\noindent{\bf Proof of Thm.\ref{zfc-t}:} By \cite[Lemma 3.1,
p.35]{NDis}, or by \cite[Lemma 2.2, Remark 2.4, p.25, p.30]{NPrep}
there is a recursive function $f:\Fm_\omega^2\to\Fm_3^2$ for which
\begin{align}\label{preserving}
\pi\models f(\f)\leftrightarrow\f,\quad f(\lnot\f)=\lnot f(\f),
\quad f(\f\lor\psi)=f(\f)\lor f(\psi),
\end{align}
for all sentences $\f,\psi$ of $\L_\omega$. Take such an $f'$ and
extend it to $\Fm_\omega$ by letting $f(\f)\de f'(\f')$ where $\f'$
is the universal closure of $\f$ (if $n$ is the smallest number such
that the free variables of $\f$ are among $v_0,...,v_n$, then $\f'$
is $\forall v_0\dots\forall v_n\f$). Then $f$ has the same
properties as $f'$ and it is defined on the whole of $\Fm_\omega$,
not only on $\Fm_\omega^2$. Define $\Tr:\Fm_\omega\to\Fmd_3$ by
$$ \Tr(\f)\de \h(f(\f))$$
for all $\f\in\Fm_\omega$. We show that this translation function
satisfies the requirements of Thm.\ref{zfc-t}. First, $\Tr$ is
recursive because both $f$ and $\h$ are such. We defined the
parameters $p_0, p_1$ so that
\begin{align}\label{pi}
\ZF\models\pi^+
\end{align} holds. Thus
$\ZF\models\ZF+\pi^+\models \f\leftrightarrow f(\f)$ by the chosen
properties of $f$, and then
$\ZF+\Delta\models\ZF+\pi^+\land\Delta\models \f\leftrightarrow \h
f\f$ by Thm.\ref{h-t}(ii), so
$\ZF+\Delta\models\f\leftrightarrow\Tr\f$ for all sentences $\f$ of
$L_\omega$. This is Thm.\ref{zfc-t}(ii). To prove
Thm.\ref{zfc-t}(i), first we show
\begin{align}\label{f}
\ZF\models\f\quad\Leftrightarrow\quad f(\ZF)+\pi^+\models f(\f).
\end{align}
Indeed, assume $\ZF\models\f$ and let $\Mm\models f(\ZF)+\pi^+$.
Then $\Mm\models\ZF$ by (\ref{preserving}), so $\Mm\models\f+\pi^+$
by our assumption $\ZF\models\f$, so $\Mm\models f(\f)$ by
(\ref{preserving}). This shows $f(\ZF)+\pi^+\models f(\f)$.
Conversely, assume now the latter, and we want to prove
$\ZF\models\f$. Let $\Mm\models\ZF$, then $\Mm\models\ZF+\pi^+$ by
(\ref{pi}), thus $\Mm\models f(\ZF)+\pi^+$ by (\ref{preserving}),
then $\Mm\models f(\f)+\pi^+$ by our assumption, and then
$\Mm\models \f$ by (\ref{preserving}) again. We have shown that
(\ref{f}) holds.

By combining this (\ref{f}) with Thm.\ref{h-t}(i) we get
$\ZF\models\f$ iff $\Tr(\ZF)\vd\Tr\f$ which is Thm.\ref{zfc-t}(i).

Later, in section \ref{free-sec}, we will also need the following:
\begin{align}\label{bool}
\Tr\mbox{ is Boolean preserving}\quad\mbox{ and }\quad
\SAx^*\vd\Tr(\lnot\f)\to\lnot\Tr\f .
\end{align}
Indeed, this follows from Thm.\ref{h-t}(iii) and from
(\ref{preserving}): Let $\f,\psi\in\Fm_\omega^0$. Then
$\vd\Tr(\f\land\psi)$,\ \  iff by the definition of $\Tr$\\
$\vd\h f(\f\land\psi)$,\ \  iff by (\ref{preserving})\\
$\vd\h(f\f\land f\psi)$,\ \  iff by Thm.\ref{h-t}(iii)\\
$\vd\h f\f\land\h f\psi$,\ \ iff by the definition of $\Tr$\\
$\vd\Tr\f\land\Tr\psi$.

Similarly,\\
$\vd\Tr(\f\to\psi)$,\ \  implies by the definition of $\Tr$\\
$\vd\h f(\f\to\psi)$,\ \  implies by (\ref{preserving})\\
$\vd\h(f\f\to f\psi)$,\ \  implies by Thm.\ref{h-t}(iii)\\
$\vd\h f\f\to\h f\psi$,\ \ implies by the definition of $\Tr$\\
$\vd\Tr\f\to\Tr\psi$.

Finally, $\SAx^*\vd\h(\lnot f\f)\to\lnot\h(f\f)$,  by
Thm.\ref{h-t}(iii), so $\SAx^*\vd\h(f(\lnot \f))\to\lnot\h(f\f)$ by
(\ref{preserving}), and then $\SAx^*\vd\Tr(\lnot\f)\to\lnot\Tr(\f)$,
as was to be shown. \qquad{\bf QED(Thm.\ref{zfc-t})}\bigskip

\section{Free algebras}\label{free-sec}
In this section we prove that the one-generated free $\Df_3$ is not
atomic. In algebraic logic, in the duality between algebras and
logics, atomicity of the Lindenbaum-Tarski algebras correspond to
G\"odel incompleteness theorem, see e.g., \cite[sec 1.4]{NPrep},
\cite{Gyenis}.

\begin{thm}\label{free-t}
The one-generated free $\Df_3$ is not atomic.
\end{thm}

\noindent{\bf Proof.} It is enough to show that the zero-dimensional
part of the free $\Df_3$ is not atomic by \cite[1.10.3(i)]{HMTII}:
\begin{align}
\gZd\gFr_1\Df_3\mbox{ is not atomic implies that }\gFr_1\Df_3\mbox{
is not atomic.}
\end{align}
Let us define $\equiv$ as $\equiv_{\emptyset}$, i.e.
$$\f\equiv\psi\qquad\mbox{ iff }\qquad \vd\f\leftrightarrow\psi .$$
Let $\gFmd_3$ and $\gFmd_3^0$ denote the word-algebra of $\Fmd_3$
and the word-algebra of sentences of $\Fmd_3$, respectively (the
latter under the operations of $\lor, \lnot$). It is easy to show
that
\begin{align}
\gFr_1\Df_3\mbox{ is isomorphic to }\gFmd_3\slash\equiv \mbox{ and
}\\
\gZd\gFr_1\Df_3\mbox{ is isomorphic to }\gFmd_3^0\slash\equiv .
\end{align}
There is a non-separable formula $\lambda\in\Fm_\omega^0$ which is
consistent with $\pi^+$ (with our concrete pairing formulas), by
\cite[pp.69-71]{NDis} (or equivalently by \cite[Lemma 2.7,
p.34]{NPrep}). Define
$$\psi\de\SAx^*\land\Tr\lambda .$$
Then $\psi\in\Fmd_3^0$. We will show that there is no atom below
$\psi\slash\equiv\ $\ and the latter is nonzero in
$\Fmd_3^0\slash\equiv$. Assume the contrary, i.e., that
\begin{align}\label{atom}
\mbox{$\delta\slash\equiv$\ \  is an atom below\ \
$\psi\slash\equiv$}
\end{align}
 and we will derive a contradiction. Let
 $$T\de\{\f\in\Fm_\omega^0 : \vd\delta\to\Tr\f\} .$$
 Then $T\subseteq\Fm_\omega^0$. We will show that $T$ is recursive
 and it separates the consequences of $\lambda$ from the sentences
 refutable from $\lambda$ which contradicts the choice of $\lambda$,
 namely that $\lambda$ is inseparable. From now on let $\f\in\Fm_\omega^0$ be arbitrary.
 Now, $\delta\slash\equiv$  being an atom implies
 \begin{align}\label{complement}
 \notvd\delta\to\Tr\f\quad\Leftrightarrow\quad\vd\delta\to\lnot\Tr\f .
 \end{align}
 From (\ref{complement}) we get that both $T$ and the complement of
 $T$ are recursively enumerable, so $T$ is recursive. Next we show
 \begin{align}\label{sep1}
 \lambda\models\f\quad\mbox{implies}\quad\f\in T .
 \end{align}
 Indeed, assume that $\lambda\models\f$. Then $\models\lambda\to\f$.
 Then, in particular, $\pi^+\models\lambda\to\f$, and so
 $\vd\Tr(\lambda\to\f)$ by Thm.\ref{zfc-t}(i). Then
 $\vd\Tr(\lambda)\to\Tr\f$ by (\ref{bool}). By Modus Ponens then
 $\SAx^*+\Tr(\lambda)\vd\Tr\f$, i.e., $\psi\vd\Tr\f$. By
 (\ref{atom}) we have
 \begin{align}\label{delta}
 \vd\delta\to\psi
 \end{align}
 so we have $\vd\delta\to\Tr\f$, i.e., $\f\in T$ as was desired.
 Next we show
 \begin{align}
 \lambda\models\lnot\f\quad\mbox{implies}\quad\f\notin T .
 \end{align}
 Indeed, assume $\lambda\models\lnot\f$. Then \\
 $(\lnot\f)\in T$ by the previous case (\ref{sep1}), and this means\\
 $\vd\delta\to\Tr(\lnot\f)$. Now, by (\ref{delta}) and the definition of $\psi$ we
 have\\
 $\vd\delta\to\SAx^*$, and thus by $\vd\delta\to\Tr(\lnot\f)$ and
 (\ref{bool})\\
 $\vd\delta\to\lnot\Tr(\f)$. Thus by (\ref{complement}) we get\\
 $\notvd\delta\to\Tr\f$, by the definition of $T$ then\\
 $\f\notin T$.

 By the above we have shown that there is no atom below
 $\psi\slash\equiv$. It remains to show that the latter is nonzero.
 This follows from the fact that $\lambda\land\pi^+$ has a model.
 Let $\Mm$ be such that $\Mm\models\lambda\land\pi^+$. Expand this
 model with $P^{\Mm}$ so that $\Mm^+\models\Delta$ also. Then
 $\Mm^+\models\Tr(\lambda)$ by Thm.\ref{zfc-t}(ii), and also
 $\Mm^+\models\SAx^*$ by $\Mm^+\models\pi^+\land\Delta$. Thus
 $\Mm^+\models\psi$, and so $\Mm^-\models\psi$ where $\Mm^-$ is the
 reduct of $\Mm^+$ to the language of $\psi$.\qquad{\bf QED}\bigskip

\bigskip\bigskip\bigskip

\noindent Hajnal Andr\'eka\ \ and\ \ Istv\'an N\'emeti\\
R\'enyi Mathematical Research Institute\\
Budapest, Re\'altanoda st.\ 13-15\\
H-1053 Hungary\\
andreka@renyi.hu, nemeti@renyi.hu


\begin{thebibliography}{}

\bibitem{dens} Andr{\'e}ka, H., Givant, S. R., Mikul{\'a}s, Sz., N{\'e}meti, I., and Simon, A.,
{\it Notions of density that imply representability in algebraic
logic}. Annals of Pure and Applied Logic 91 (1998), 93-190.

\bibitem{HbPhL} Andr{\'e}ka, H., N{\'e}meti, I., and Sain, I.,
{\it Algebraic Logic}.  In: Handbook of Philosophical Logic Vol. 2,
Second Edition, Kluwer, 2001. pp.133-296.
\url{http://www.math-inst.hu/pub/algebraic-logic/handbook.pdf}

\bibitem{GKWZ} Gabbay, D. M., Kurucz, \'A., Wolter, F., and
Zakharyaschev, M., {\it Many-dimensional modal logics: theory and
applications.}, Elsevier, 2003.

\bibitem{Gyenis} Gyenis, Z., On atomicity of free algebras in certain
cylindric-like varieties. Logic Journal of the IGPL 19,1 (2011),
44-52.

\bibitem{HMTII} Henkin, L., Monk, J. D., and Tarski, A.,
{\it Cylindric Algebras Parts I and II}, North-Holland, 1985.

\bibitem{NPrep} N\'emeti, I., {\it Logic with three variables has
G\"odel's incompleteness property - thus free cylindric algebras are
not atomic.}, Mathematical Institute of the Hungarian Academy of
Sciences, Preprint No 49/1985, 1985.
\url{http://www.math-inst.hu/~nemeti/NDis/NPrep85.pdf}

\bibitem{NDis} N\'emeti, I., {\it Free algebras and decidability in
algebraic logic.}, Dissertation with the Hungarian Academy of
Sciences, Budapest, 1986. In Hungarian. English summary is
\cite{NSum}. \url{http://www.math-inst.hu/~nemeti/NDis/NDis86.pdf}

\bibitem{NSum} N\'emeti, I., {\it Free algebras and decidability
in algebraic logic. Summary in English.} 12pp.
\url{http://www.math-inst.hu/~nemeti/NDis/NSum.pdf}

\bibitem{shortnote} N\'emeti, I., {\it Formalizing set theory in weak fragments
of algebraic logic (updated in June 2011)} 4pp.
\url{http://www.math-inst.hu/~nemeti/NDis/formalizingsettheory.pdf}

\bibitem{NeSiIGPL97} N\'emeti, I., and Simon, A., {\it Relation
algebras from cylindric and polyacid algebras.} Logic Journal of the
IGPL 5,4 (1997), 575-588.

\bibitem{Say04} Sayed-Ahmed, T., {\it Tarskian Algebraic Logic.},
Journal on Relation Methods in Computer Science 1 (2004), pp.3-26.

\bibitem{Say05} Sayed-Ahmed, T., {\it Algebraic logic, where does it stand today?},
The Bulletin of Symbolic Logic 11,4 (2005), pp.465-516.

\bibitem{Simon} Simon, A., {\it Connections between quasi-projective relation algebras
and cylindric algebras.}, Algebra Universalis 56,3-4 (2007),
263-301.

\bibitem{T53a} Tarski, A., {\it A formalization of set theory
without variables.}, J. Symbolic Logic 18 (1953), p.189.

\bibitem{T53b} Tarski, A., {\it An undecidable system of sentential calculus.},
J. Symbolic Logic 18 (1953), p.189.

\bibitem{TG} Tarski, A., and Givant, S. R., {\it Formalizing set
theory without variables.}, AMS Colloquium Publications Vol. 41,
AMS, Providence, Rhode Island, 1987.

\end{thebibliography}
\end{document}